\documentclass[a4paper, 11pt]{amsart}
\usepackage[utf8]{inputenc}
\usepackage[T1]{fontenc}
\usepackage{amsmath,amssymb,amsthm,amsfonts}
\usepackage[a4paper,hmargin=2.5cm,vmargin=2.5cm]{geometry}
\usepackage[english]{babel}
\usepackage[unicode]{hyperref}
\usepackage{enumitem}
\usepackage[square,numbers]{natbib}

\theoremstyle{plain}
\newtheorem{theorem}{Theorem}
\newtheorem{lemma}[theorem]{Lemma}

\newtheorem{phrase}[theorem]{Proposition}

\theoremstyle{definition}
\newtheorem{remark}[theorem]{Remark}
\newtheorem{define}[theorem]{Definition}

\newtheorem{example}[theorem]{Example}

\theoremstyle{remark}

\title{A sharp lower bound for the first nonzero Steklov eigenvalue of the Witten Laplacian on weighted surfaces with nonnegative Bakry--Emery--Ricci tensor}
\author{Artem M. Danilov
}
\date{}

\begin{document}

\maketitle

\begin{abstract}
We consider compact connected weighted surfaces with smooth boundary whose Bakry--Emery--Ricci tensor is nonnegative. We assume that the geodesic curvature of the boundary is bounded below in terms of the
normal and tangential derivatives of the weight function. Under these assumptions, we prove a lower bound for the first nonzero Steklov eigenvalue of the Witten Laplacian, which is sharp in the sense that equality holds if and only if the surface is isometric to a Euclidean disk and the weight function is constant. 
\end{abstract}

\section{Introduction}

The main goal of this paper is to study the first nonzero Steklov eigenvalue of the Witten Laplacian on a compact connected Riemannian surface with smooth boundary endowed with a measure. Throughout this paper, by a \emph{surface} we mean a two-dimensional smooth manifold, possibly with boundary.

\begin{define}[see, e.g., \cite{item1}]
Let $(M, g)$ be a compact connected smooth Riemannian manifold endowed with a measure $\sigma\, dv_g$, where $\sigma = e^{-f} \in C^{\infty}(M)$ is a positive density and $dv_g$ is the Riemannian measure induced by the metric $g$. Such a manifold is called \emph{weighted} and is denoted by $(M, g, \sigma)$. The \emph{Witten Laplacian} (also called the Bakry--Emery Laplacian or the weighted Laplacian) is defined by
\[
L_{\sigma} = \Delta_g - \frac{1}{\sigma}g(\nabla^g \sigma, \nabla^g ) = \Delta_g + g(\nabla^g f , \nabla^g),
\]
where $\Delta_g = -\mathrm{div} \, \nabla^g$ is the standard Laplace--Beltrami operator and $\nabla^g$ is the gradient with respect to the metric $g$.
\end{define}
The Steklov problem for the Witten Laplacian is the eigenvalue problem
\begin{equation}\label{eq:1.1}
\begin{cases}
L_{\sigma} u = 0 & \text{in } M, \\
\dfrac{\partial u}{\partial \vec{n}} = \lambda u & \text{on } \partial M,
\end{cases}
\end{equation}
where $\vec{n}$ is the outward unit normal vector field. In this problem $u$ and
$\lambda$ are called an eigenfunction and an eigenvalue, respectively. 

It is known (see, e.g., \cite{jammes}) that the set of
eigenvalues of \eqref{eq:1.1} is a discrete nonnegative spectrum which, counted with
multiplicity, forms a sequence
\[
0 = \lambda_0 < \lambda_1\leq\lambda_2\leq\lambda_3 \leq \ldots \;.
\]
This is the spectrum of the Dirichlet-to-Neumann operator
$D_{\sigma}\colon H^{\frac{1}{2}}(\partial M) \to H^{-\frac{1}{2}}(\partial M)$,
$D_{\sigma}(h) = \tfrac{\partial P_{\sigma}h}{\partial \vec{n}}$, where $P_{\sigma}h$
is the solution of the problem
\begin{equation*}\label{eq:1.2}
\begin{cases}
L_{\sigma} v = 0 & \text{in } M, \\
v = h & \text{on } \partial M.
\end{cases}
\end{equation*}
The operator
$D_\sigma$ is self-adjoint with respect to the inner product of $L^2(\partial M, e^{-f}\,ds)$. Moreover, the spectrum has a variational characterization. In particular,
\[
\lambda_1=\inf_{\substack{u\in H^1(M)\\ \int_{\partial M}u\,e^{-f}dv_g=0}}
\frac{\int_M g(\nabla^g u, \nabla^g u)\,e^{-f}\,dv_g}{\int_{\partial M}u^2e^{-f}\,ds}.
\]
where $ds$ is the Riemannian measure on $\partial M$ induced by $g$.

There are several known estimates for the first nonzero Steklov eigenvalue.
For a simply connected domain $\Omega$ in $\mathbb{R}^2$ with $C^2$ boundary,
Payne~\cite{payne}, using a result of Weinstock~\cite{weinstock}, obtained
sharp lower and upper bounds for the first nonzero Steklov eigenvalue in terms
of the geodesic curvature of the boundary.

Since we work with smooth boundaries throughout, we state Payne's result
in this setting only.

\begin{theorem}[{\cite[Theorem 1]{payne}}]\label{payne_theor}
Let $\Omega$ be a simply connected bounded domain in $\mathbb{R}^2$ with smooth boundary $\partial \Omega$. Then the first nonzero Steklov eigenvalue $\lambda_1$ satisfies
\begin{equation*}
\max_{x \in \partial \Omega} k_g(x) \geq \lambda_1 \geq \min_{x \in \partial \Omega} k_g(x),
\end{equation*}
where $k_g$ is the geodesic curvature of the boundary. Moreover, both equalities hold if and only if $\Omega$ is a Euclidean disk.
\end{theorem}

Later, Escobar \cite{escobar} generalized the lower bound of Theorem~\ref{payne_theor} to Riemannian surfaces.

\begin{theorem}[Escobar~{\cite[Theorem 1]{escobar}}]\label{escobar_theor}
Let $(M, g)$ be a compact surface with smooth boundary $\partial M$. Assume that $M$ has nonnegative Gaussian curvature $K$ and that the geodesic curvature $k_g$ of $\partial M$ satisfies $k_g \geq c > 0$ for some constant $c$. Then the first nonzero Steklov eigenvalue $\lambda_1$ satisfies
\begin{equation*}
\lambda_1 \geq c.
\end{equation*}
Moreover, equality holds if and only if $M$ is isometric to a Euclidean disk of radius $c^{-1}$.
\end{theorem}

The main result of this paper is Theorem \ref{main_t} in Subsection~\ref{main_theorem}, which is a weighted version of Theorem~\ref{escobar_theor}. The proof is based on the Bochner formula for the Bakry--Emery--Ricci tensor and the strong maximum principle, and it generalizes the computations used by Cho and Seo in~\cite{cho}.

In~\cite{cho} the problem
\begin{equation*}
\begin{cases}
H_\sigma v=0 & \text{in } M,\\
\dfrac{\partial v}{\partial \vec{n}}=\lambda v & \text{on } \partial M,
\end{cases}
\end{equation*}
where $H_\sigma=\Delta_g+\tfrac14|\nabla^g f|^2+\tfrac12\Delta_g f$ is the Schr\"odinger operator, is considered as a particular case. A problem involving this Schr\"odinger operator can also be obtained from problem~\eqref{eq:1.1}. Indeed, $L_\sigma$ is self-adjoint on $L^2(\sigma\,dv_g)$ and is unitarily equivalent, through the transform $\sqrt{\sigma}\colon L^2(\sigma\,dv_g)\to L^2(dv_g)$, to the Schr\"odinger operator $H_\sigma$. However, in the presence of a boundary, this map also changes the boundary condition. Namely, for $v=\sqrt{\sigma}\,u$ the Steklov condition $\partial u/\partial\vec n=\lambda u$ turns into 
\begin{equation*}
\frac{\partial v}{\partial \vec{n}}+\frac12\,\frac{\partial f}{\partial \vec{n}}\,v=\lambda v
\quad \text{on } \partial M,
\end{equation*}
which is not the Steklov condition considered in~\cite{cho}. Thus Theorem~\ref{main_t} does not follow from the results of~\cite{cho}.
\subsection*{Acknowledgements}
The author is grateful to his scientific advisor, Aleksei V.~Penskoi,
for valuable discussions, helpful remarks on the text, and constant
attention to this work.
\section{Preliminaries}

\begin{define}[see, e.g., {\cite[Introduction]{ma}}]
Let $(M, g, \sigma)$ be a compact connected weighted surface and let $\sigma = e^{-f} \in C^{\infty}(M)$. The \emph{Bakry--Emery--Ricci tensor} is defined by $\mathrm{Ric}_{\sigma} = \mathrm{Ric} + \mathrm{Hess}\,f$.
\end{define}

Since $\partial M$ is a one-dimensional smooth compact manifold without boundary,
the classification theorem for one-dimensional manifolds (see, e.g.,
\cite[Appendix]{milnor}) implies that each connected component of $\partial M$
is diffeomorphic to $S^1$. Compactness of $\partial M$ ensures that the number
of such components is finite. From this we deduce the following.

\begin{phrase}\label{boundary_cor}
Let $(M, g, \sigma)$ be a compact connected weighted surface with a smooth boundary, and let $\sigma = e^{-f} \in C^{\infty}(M)$. Then each
connected component $\Sigma \subset \partial M$ is a smooth simple closed curve. In
particular, each component admits an arc-length parametrization
\[
\gamma_{\Sigma} \colon [0, \ell] \to \Sigma, \quad |\dot{\gamma_{\Sigma}}(t)|_g = 1
\quad \forall\, t \in [0, \ell],
\]
where $\ell$ is the length of the curve $\gamma_{\Sigma}$, with $\gamma_{\Sigma}(0)
= \gamma_{\Sigma}(\ell)$ and $\dot{\gamma_{\Sigma}}(0) = \dot{\gamma_{\Sigma}}(\ell)$.
The velocity vectors $\dot{\gamma_{\Sigma}}(t)$ define a smooth unit tangent
vector field along $\partial M$. Moreover, any smooth unit tangent vector field
$\vec{\nu}$ along $\partial M$ coincides with $\pm\dot{\gamma_{\Sigma}}$ on each
component $\Sigma$. Hence $\left|\dfrac{\partial f}{\partial \vec{\nu}}(x)\right|$
does not depend on the choice of $\vec{\nu}$.
\end{phrase}
\begin{define}
For $a,b\in\mathbb{R}$, set
\[
\mu_+(a,b)=\frac{a+\sqrt{a^2+b^2}}{2}.
\]
\end{define}
\begin{remark}\label{mu}
Note that $\mu_+\geq0$ and $\mu_+\geq a$.
Since $\mu_+$ depends on $b$ only through $b^2$, the number
$\mu_+\!\left(\dfrac{\partial f}{\partial \vec{n}},\dfrac{\partial f}{\partial \vec{\nu}}\right)$ does not
depend on the choice of $\vec\nu$ in Proposition~\ref{boundary_cor}.
\end{remark}
\section{Main results}

\subsection{Main theorem}
\label{main_theorem}

\begin{theorem}\label{main_t}
Let $(M, g, \sigma)$ be a compact connected weighted surface with a smooth boundary, and let $\sigma = e^{-f} \in C^{\infty}(M)$. Assume that
\begin{equation}\label{eq:1.3}
\mathrm{Ric}_{\sigma, x}(X,X)\geq 0 \quad \forall x \in M,\ \forall X \in T_xM.
\end{equation}
Assume also that the geodesic curvature $k_g$ of $\partial M$ satisfies
\begin{equation}\label{eq:1.4}
k_g(x) - \mu_+\!\left(\dfrac{\partial f}{\partial \vec{n}}(x),\dfrac{\partial f}{\partial \vec{\nu}}(x)\right) \geq c > 0 \quad \forall x \in \partial M,
\end{equation}
where $\vec{n}$ is the outward unit normal vector field, $\vec{\nu}$ is a unit tangent vector field along $\partial M$, and $c$ is a constant. Then the first nonzero eigenvalue $\lambda_1$ of the Steklov problem~\eqref{eq:1.1} satisfies $\lambda_1 \geq c$. Moreover, equality holds if and only if $M$ is isometric to a Euclidean disk of radius $c^{-1}$ and $f$ is constant.
\end{theorem}
\begin{example}
Let $M$ be the unit Euclidean disk in $\mathbb{R}^2$ and $f(x,y)=\varepsilon x$ with $0<\varepsilon<1$. Then $\mathrm{Ric}_\sigma=0$, $k_g=1$, and $\mu_+\le\varepsilon$ on $\partial M$, so assumptions~\eqref{eq:1.3} and~\eqref{eq:1.4} hold with $c=1-\varepsilon$, and $\lambda_1\geq1-\varepsilon$.
\end{example}
Condition~\eqref{eq:1.4} may look involved. However, the estimate is sharp, since among all weighted surfaces satisfying~\eqref{eq:1.3} and~\eqref{eq:1.4}, equality $\lambda_1=c$ is attained only by those isometric to a Euclidean disk of radius $c^{-1}$ with a constant weight function $f$. This is the most natural geometric case. We do not claim that condition~\eqref{eq:1.4} is optimal for a nonconstant weight. Note also that for $f \equiv \mathrm{const}$, condition~\eqref{eq:1.4} reduces to $k_g \geq c > 0$, and Theorem~\ref{main_t} recovers the classical result of Escobar~\cite{escobar}.
\subsection{Smoothness of eigenfunctions}
A priori, an eigenfunction of~\eqref{eq:1.1} is only assumed to lie in $H^1(M)$, so the problem is understood in the weak sense.

\begin{define}\label{def7}
A function $u\in H^1(M)$ is called a \emph{weak eigenfunction} of the Steklov problem~\eqref{eq:1.1} with eigenvalue $\lambda$ if
\begin{equation}\label{eq:weak_steklov}
\int_M \sigma\, g(\nabla^g u,\nabla^g\psi)\,dv_g=\lambda\int_{\partial M}\sigma\, u\,\psi\,ds\qquad\forall\psi\in H^1(M).
\end{equation}
\end{define}

\begin{define}\label{def8}
A function $u\in C^2(M)$ is called a \emph{classical eigenfunction} of~\eqref{eq:1.1} with eigenvalue $\lambda$ if it satisfies the following system of equations pointwise
\begin{equation}\label{eq:classical_steklov}
\begin{cases}
L_\sigma u=0 & \text{in } M, \\
\dfrac{\partial u}{\partial\vec n}=\lambda u & \text{on } \partial M.
\end{cases}
\end{equation}
\end{define}

We prove below that every weak eigenfunction is smooth and is a classical eigenfunction. Conversely, by Green's identity, every classical smooth eigenfunction is a weak one. Thus the two notions coincide, and in the rest of the paper we may use classical computations.

Taking $\psi=e^{f}\varphi$ in Definition~\ref{def7} for $\varphi\in H^1(M)$, we obtain the following equivalent definition of a weak solution of the Steklov problem~\eqref{eq:1.1}.

\begin{define}\label{def9}
A function $u\in H^1(M)$ is a \emph{weak eigenfunction} of~\eqref{eq:1.1} with eigenvalue $\lambda$ if
\begin{equation}\label{eq:weak_neumann}
\int_M g(\nabla^g u,\nabla^g\varphi)\,dv_g=-\int_M g(\nabla^g f,\nabla^g u)\varphi\,dv_g+\lambda\int_{\partial M}u\varphi\,ds\qquad\forall\varphi\in H^1(M).
\end{equation}
\end{define}

\begin{remark}\label{rem:compat}
Taking $\varphi\equiv1$ in~\eqref{eq:weak_neumann}, we see that every weak eigenfunction satisfies
\[
\int_M g(\nabla^g f,\nabla^g u)\,dv_g=\lambda\int_{\partial M}u\,ds .
\]
\end{remark}

\begin{phrase}[{\cite[Lemma~28.4.1]{Mendoza}}]\label{neumann_reg}
Let $(M,g)$ be a smooth compact connected Riemannian manifold with boundary. For $a\geq 0$, let $F\in H^a(M)$ and $h\in H^{a+1/2}(\partial M)$ satisfy $\int_M F\,dv_g=\int_{\partial M}h\,ds$. Then there exists a unique mean-zero function $q\in H^{a+2}(M)$ such that
\[
\begin{cases}
\Delta_g q=-F & \text{in } M, \\
\dfrac{\partial q}{\partial\vec n}=h & \text{on } \partial M.
\end{cases}
\]
\end{phrase}

The following Sobolev embedding is a direct consequence of Theorem~2.30 in~\cite{Aubin}.

\begin{phrase}\label{sobolev_emb}
Let $N$ be a smooth compact manifold of dimension $n$ with smooth boundary. Then for every positive integer $k$ such that $k+\frac{n}{2}<s$, there is a continuous embedding
\[
H^s(N) \subset C^k(N).
\]
\end{phrase}

\begin{phrase}\label{smoothness}
Let $(M, g, \sigma)$ be a compact connected weighted surface with a smooth boundary, let $\sigma = e^{-f} \in C^{\infty}(M)$, and let $u \in H^1(M)$ be a weak eigenfunction of the Steklov problem~\eqref{eq:1.1}. Then $u \in C^{\infty}(M)$, and $u$ is a classical eigenfunction.
\end{phrase}
\begin{proof}
Put $F:=g(\nabla^g f,\nabla^g u)$ and $h:=\lambda\, u|_{\partial M}$. Then Definition~\ref{def9} gives
\begin{equation}\label{eq:weak_Fh}
\int_M g(\nabla^g u,\nabla^g\varphi)\,dv_g=-\int_M F\varphi\,dv_g+\int_{\partial M}h\varphi\,ds\qquad\forall\varphi\in H^1(M),
\end{equation}
and by Remark~\ref{rem:compat} we get $\int_M F\,dv_g=\int_{\partial M}h\,ds$.

We prove by induction on $k\geq 0$ that $u\in H^{k+1}(M)$. The case $k=0$ is the assumption. Let $u\in H^{k+1}(M)$. Since $f$ is smooth, $F\in H^k(M)$, and $h\in H^{k+1/2}(\partial M)$ by the trace theorem (see, e.g., \cite[Theorem~5.1.7]{Agranovich_Sobolev}). Hence, by Proposition~\ref{neumann_reg} with $a=k$, there exists $q\in H^{k+2}(M)$ with $\Delta_g q=-F$ in $M$ and $\dfrac{\partial q}{\partial\vec n}=h$ on $\partial M$. By Green's identity, which extends from smooth functions to $H^2(M)$ by density, $q$ also satisfies~\eqref{eq:weak_Fh}.

Subtracting~\eqref{eq:weak_Fh} for $q$ from~\eqref{eq:weak_Fh} for $u$, we get, for $w:=u-q\in H^1(M)$,
\[
\int_M g(\nabla^g w,\nabla^g\varphi)\,dv_g=0\qquad\forall\varphi\in H^1(M).
\]
Hence $w$ is constant since $\nabla^g w=0$ and $M$ is connected. Therefore $u=q+c'\in H^{k+2}(M)$.

Thus $u\in H^m(M)$ for every $m\geq1$, and by Proposition~\ref{sobolev_emb} with $n=2$ we get $u\in C^\infty(M)$.

It remains to show that $u$ is a classical eigenfunction. Since $u$ is smooth, Green's identity gives, for every $\varphi\in C^\infty(M)$,
\[
\int_M g(\nabla^g u,\nabla^g\varphi)\,dv_g=\int_M(\Delta_g u)\,\varphi\,dv_g+\int_{\partial M}\dfrac{\partial u}{\partial\vec n}\,\varphi\,ds .
\]
Substituting this into~\eqref{eq:weak_Fh}, we obtain
\[
\int_M(\Delta_g u+F)\,\varphi\,dv_g+\int_{\partial M}\Bigl(\dfrac{\partial u}{\partial\vec n}-h\Bigr)\varphi\,ds=0\qquad\forall\varphi\in C^\infty(M).
\]
For $\varphi\in C_c^\infty(\mathrm{Int}\,M)$, the boundary integral vanishes, and since $\Delta_g u+F$ is continuous, we get $\Delta_g u=-F$ pointwise in $M$, i.e.,\ $L_\sigma u=0$. Then the first integral vanishes for all $\varphi\in C^\infty(M)$, so
\[
\int_{\partial M}\Bigl(\dfrac{\partial u}{\partial\vec n}-h\Bigr)\varphi\,ds=0\qquad\forall\varphi\in C^\infty(M).
\]
Since $\varphi|_{\partial M}$ can be an arbitrary smooth function on $\partial M$ and $\dfrac{\partial u}{\partial\vec n}-h$ is continuous, we get $\dfrac{\partial u}{\partial\vec n}=h=\lambda u$ pointwise on $\partial M$. Thus $u$ is a classical eigenfunction.
\end{proof}
\subsection{Proof of the main theorem}

The proof relies on the Bochner formula for the Bakry--Emery--Ricci tensor and the strong maximum principle for the Witten Laplacian.

\begin{phrase}[Bochner formula for the Bakry--Emery--Ricci tensor, {\cite[Lemma~2.1]{gonzalez}}]\label{bochner}
Let $(M, g, \sigma)$ be a compact connected weighted surface with a smooth boundary, and let $\sigma = e^{-f} \in C^{\infty}(M)$. Then for every function $\varphi \in C^{\infty}(M)$ one has
\[
\frac{1}{2} L_{\sigma}|\nabla^g \varphi|^2 = -|\mathrm{Hess}\,\varphi|^2+g(\nabla^g \varphi, \nabla^g (L_{\sigma}\varphi))-\mathrm{Ric}_{\sigma}(\nabla^g \varphi, \nabla^g \varphi),
\]
where $|X|^2 := g(X, X)$ for any tangent vector $X$.
\end{phrase}

As a direct consequence of Theorem~$3.5$ and Lemma~$3.4$ in~\cite{gilbarg} applied to the Witten Laplacian, we get the following statement.

\begin{phrase}[Strong maximum principle]\label{strong_max}
Let $(M, g, \sigma)$ be a compact connected weighted surface with a smooth boundary, and let $\sigma = e^{-f} \in C^{\infty}(M)$. Let $\varphi \in C^{\infty}(M)$ satisfy $L_{\sigma}\varphi \leq 0$ in $\mathrm{Int}\, M$. Then either $\varphi$ is constant on $M$, or $\varphi$ attains its maximum
only at points of $\partial M$, and at every maximum point
$p \in \partial M$ one has $\dfrac{\partial \varphi}{\partial \vec{n}}(p) > 0$,
where $\vec{n}$ is the outward unit normal vector field.
\end{phrase}

Let $(M, g, \sigma)$ satisfy the assumptions of Theorem~\ref{main_t}, and let $u$ be an eigenfunction of the Steklov problem~\eqref{eq:1.1} corresponding to the first nonzero eigenvalue $\lambda_1$. By Proposition~\ref{bochner},
\begin{equation}\label{eq:1.0}
\begin{gathered}
\frac{1}{2} L_{\sigma}|\nabla^g u|^2 = -|\mathrm{Hess}\,u|^2-\mathrm{Ric}(\nabla^g u, \nabla^g u)-\mathrm{Hess}\,f(\nabla^g u, \nabla^g u)\\=-|\mathrm{Hess}\,u|^2-\mathrm{Ric}_{\sigma}(\nabla^g u, \nabla^g u) \leq0,
\end{gathered}
\end{equation}
where the last inequality follows from assumption~\eqref{eq:1.3} of Theorem~\ref{main_t}. Hence $L_{\sigma}|\nabla^g u|^2 \leq 0$.

By Proposition~\ref{smoothness}, $u \in C^{\infty}(M)$, and therefore $|\nabla^g u|^2\in C^{\infty}(M)$. Moreover, since $M$ is compact, the continuous function $|\nabla^g u|^2$ attains its maximum at some point $p$. By Proposition~\ref{strong_max}, two cases are possible:
\begin{enumerate}[label=\Roman*., ref=\Roman*]
  \item \label{I} $p \in \partial M$, the function $|\nabla^g u|^2$ is nonconstant on $M$, and $\dfrac{\partial |\nabla^g u|^2}{\partial \vec{n}}(p)>0$;
  \item \label{II} $|\nabla^g u|^2$ is constant on $M$.
\end{enumerate}

We first prove the inequality in Case~\ref{I}. We choose a local orthonormal frame $\{e_1, e_2\}$ in a neighborhood of $p$ such that $e_1$ is tangent to $\partial M$ and $e_2 = \vec{n}$ is the outward normal. Note that $(e_2f)(x) = \dfrac{\partial f}{\partial \vec{n}}(x)$, and Proposition~\ref{boundary_cor} implies that $\left\lvert (e_1f)(x)\right\rvert = \left\lvert \dfrac{\partial f}{\partial \vec{\nu}}(x)\right\rvert$ for every point $x \in \partial M$. Note that the geodesic curvature $k_g$ of $\partial M$ satisfies
\begin{equation}\label{eq:1.5}
\nabla_{e_1} e_1 = -k_g\, e_2,
\end{equation}
where $\nabla$ is the Levi-Civita connection on $M$. This immediately implies
\begin{equation}\label{eq:1.6}
\nabla_{e_1} e_2 = k_g\, e_1,
\end{equation}
because
\[
g(\nabla_{e_1} e_2, e_1) = e_1\bigl(g(e_2, e_1)\bigr) - g(e_2, \nabla_{e_1} e_1) = k_g,
\]
and
\[
g(\nabla_{e_1} e_2, e_2) = \frac{1}{2}\,e_1\bigl(g(e_2, e_2)\bigr) = 0.
\]

Since $p$ is a maximum point of the function $|\nabla^g u|^2$, its tangential derivative vanishes at $p$. Consequently, at $p \in \partial M$, we have
\begin{equation}\label{eq:1.7}
\begin{gathered}
0 = e_1\bigl(|\nabla^g u|^2\bigr)
= 2g(\nabla_{e_1}\nabla^g u, \nabla^g u)
= 2\,\mathrm{Hess}\,u(e_1, \nabla^g u)  \\
= 2\bigl(g(\nabla^g u, e_1)\,\mathrm{Hess}\,u(e_1,e_1)
+ g(\nabla^g u, e_2)\,\mathrm{Hess}\,u(e_1,e_2)\bigr).
\end{gathered}
\end{equation}
Using \eqref{eq:1.6} and the fact that $u$ is an eigenfunction of the Steklov problem~\eqref{eq:1.1}, we obtain
\begin{gather}\label{eq:1.8}
\mathrm{Hess}\,u(e_1,e_2)
  = e_1 (e_2 u) - (\nabla_{e_1} e_2)\,u
  = e_1(\lambda_1 u) - k_g\,e_1 u
  = (\lambda_1 - k_g)\,g(\nabla^g u, e_1).
\end{gather}
Combining \eqref{eq:1.7} and \eqref{eq:1.8}, we get
\[
  g(\nabla^g u, e_1)
  \bigl(\mathrm{Hess}\,u(e_1,e_1) + (\lambda_1 - k_g)\,g(\nabla^g u, e_2)\bigr) = 0,
\]
which splits into two cases:
\begin{enumerate}[label=(\roman*)]
  \item \label{i} $g(\nabla^g u, e_1)(p) = 0$;
  \item \label{ii}$\mathrm{Hess}\,u(e_1,e_1)(p) = (k_g(p) - \lambda_1) g(\nabla^g u, e_2)(p)$.
\end{enumerate}

Consider the case \ref{i}. It follows from \eqref{eq:1.8} that $\mathrm{Hess}\,u(e_1,e_2) = 0$ at $p$.
Since $p$ is a maximum point of the function $|\nabla^g u|^2$, we have
\[
  e_1 (e_1(|\nabla^g u|^2)) \leq 0 \quad \text{at } p.
\]
Expanding the left-hand side, we obtain, at $p \in \partial M$,
\begin{equation}\label{eq:1.9}
\begin{gathered}
0 \geq  \frac{1}{2}e_1 (e_1(|\nabla^g u|^2))
= e_1\bigl(\mathrm{Hess}\,u(e_1,\nabla^g u)\bigr) \\
= e_1\bigl(g(\nabla^g u,e_1)\,\mathrm{Hess}\,u(e_1,e_1)
   + g(\nabla^g u,e_2)\,\mathrm{Hess}\,u(e_1,e_2)\bigr)
    \\
= \mathrm{Hess}\,u(e_1,e_1)^2
   + g(\nabla^g u,\nabla_{e_1}e_1)\,\mathrm{Hess}\,u(e_1,e_1)
 + g(\nabla^g u,e_2)\,e_1\bigl(\mathrm{Hess}\,u(e_1,e_2)\bigr)
   ,
\end{gathered}
\end{equation}
where we used the equalities $g(\nabla^g u,e_1) = 0$ and $\mathrm{Hess}\,u(e_1,e_2) = 0$ at $p$.

Since $u$, $k_g$, $e_1$, and $e_2$ are smooth near $p$, and $e_1=\pm\dot\gamma$ for the arc-length parametrization $\gamma$ of $\partial M$ by Proposition~\ref{boundary_cor}, we can differentiate \eqref{eq:1.8} along $\gamma$, i.e., in the direction $e_1$, at $p\in\partial M$. We obtain
\begin{equation}\label{eq:1.11}
\begin{gathered}
  e_1\bigl(\mathrm{Hess}\,u(e_1,e_2)\bigr)
  = (\lambda_1 - k_g)\,\bigl(\mathrm{Hess}\,u(e_1,e_1)
    + g(\nabla^g u,\nabla_{e_1}e_1)\bigr)-e_1(k_g)\, g(\nabla^g u, e_1)  \\
  = (\lambda_1 - k_g)\,\bigl(\mathrm{Hess}\,u(e_1,e_1)
    + g(\nabla^g u,\nabla_{e_1}e_1)\bigr),
\end{gathered}
\end{equation}
where we used that $\lambda_1$ is constant and that $g(\nabla^g u,e_1)(p)=0$.

Substituting \eqref{eq:1.11} into \eqref{eq:1.9} and applying \eqref{eq:1.5}, we obtain the inequality
\begin{equation}\label{eq:1.12}
\begin{gathered}
  0 \geq \mathrm{Hess}\,u(e_1,e_1)^2- k_g\,g(\nabla^g u,e_2)\,\mathrm{Hess}\,u(e_1,e_1)
\\ +  g(\nabla^g u,e_2)\,(\lambda_1 - k_g)\,\bigl(\mathrm{Hess}\,u(e_1,e_1)
           - k_g\, g(\nabla^g u,e_2)\bigr)
\end{gathered}
\end{equation}
at $p \in \partial M$.

On the other hand, since $\dfrac{\partial |\nabla^g u|^2}{\partial \vec{n}}(p)>0$ and $g(\nabla^g u, e_1)(p) = 0$, we have
\begin{equation}\label{eq:1.13}
\begin{gathered}
  0 < \frac{1}{2}\,e_2\bigl(g(\nabla^g u,\nabla^g u)\bigr)
     = \mathrm{Hess}\,u(\nabla^g u,e_2)\\
     = g(\nabla^g u,e_2)\,\mathrm{Hess}\,u(e_2,e_2) +  g(\nabla^g u,e_1)\,\mathrm{Hess}\,u(e_1,e_2)\\
    = g(\nabla^g u,e_2)\,\bigl(-\Delta_g u - \mathrm{Hess}\,u(e_1,e_1)\bigr) \\= g(\nabla^g u,e_2)\,\bigl(g(\nabla^g f, \nabla^g u) - \mathrm{Hess}\,u(e_1,e_1)\bigr)
\end{gathered}
\end{equation}
at $p \in \partial M$, where the last equality follows from the fact that $u$ is an eigenfunction of the Steklov problem~\eqref{eq:1.1}.

Applying \eqref{eq:1.12} and \eqref{eq:1.13}, we conclude that
\begin{equation}\label{eq:1.14}
\begin{gathered}
  \mathrm{Hess}\,u(e_1,e_1)^2
  + 2(\lambda_1 - k_g)\,g(\nabla^g u,e_2)\,\mathrm{Hess}\,u(e_1,e_1)
   - k_g(\lambda_1 - k_g)\,g(\nabla^g u,e_2)^2
  \\
  \leq\lambda_1\,g(\nabla^g u,e_2)\,\mathrm{Hess}\,u(e_1,e_1)<
 \lambda_1\,g(\nabla^g f, \nabla^g u)\, g(\nabla^g u,e_2)
\end{gathered}
\end{equation}
at $p \in \partial M$. Finally, inequality~\eqref{eq:1.14}, the condition $g(\nabla^g u, e_1)(p) = 0$, and the fact that $u$ is an eigenfunction of the Steklov problem~\eqref{eq:1.1} give
\begin{equation*}
\begin{gathered}
  0 \leq \left(\mathrm{Hess}\,u(e_1,e_1) + (\lambda_1 - k_g)\, g(\nabla^g u,e_2)\right)^2 \\
    < (\lambda_1 - k_g)^2\, g(\nabla^g u,e_2)^2
       + k_g(\lambda_1 - k_g)\, g(\nabla^g u,e_2)^2
       +\lambda_1\,g(\nabla^g f, \nabla^g u)\, g(\nabla^g u,e_2)\\
    = g(\nabla^g u,e_2)^2\,\bigl(\lambda_1(\lambda_1 - k_g)
       + \lambda_1 \, g(\nabla^g f,e_2)\bigr)
    = g(\nabla^g u,e_2)^2 \lambda_1\!\left(\lambda_1 - k_g + \dfrac{\partial f}{\partial \vec{n}}\right)
\end{gathered}
\end{equation*}
at $p \in \partial M$. Therefore, from the assumption \eqref{eq:1.4}, we get the inequality
\begin{equation*}
\begin{gathered}
\lambda_1>k_g(p)-\dfrac{\partial f}{\partial\vec n}(p)\geq k_g(p)-\mu_+\!\left(\dfrac{\partial f}{\partial\vec n}(p),\dfrac{\partial f}{\partial\vec\nu}(p)\right)\geq c,
\end{gathered}
\end{equation*}
where we used that $\mu_+(a,b)\ge a$.

Now consider case \ref{ii}. In this case, at the point $p$ we have
\begin{equation}\label{eq:1.15}
\begin{gathered}
  \mathrm{Hess}\,u(e_1,e_1) = (k_g - \lambda_1)\, g(\nabla^g u,e_2).
\end{gathered}
\end{equation}
Since $u$ is an eigenfunction of the Steklov problem~\eqref{eq:1.1}, the inequality $\dfrac{\partial |\nabla^g u|^2}{\partial \vec{n}}(p)>0$ together with \eqref{eq:1.8} and~\eqref{eq:1.15} gives
\begin{equation}\label{eq:1.16}
\begin{gathered}
  0 < \frac{1}{2}\,e_2\bigl(g(\nabla^g u, \nabla^g u)\bigr)=\mathrm{Hess}\,u(e_2,\nabla^g u) \\
    = g(\nabla^g u,e_1)\,\mathrm{Hess}\,u(e_1,e_2)
       + g(\nabla^g u,e_2)\,\mathrm{Hess}\,u(e_2,e_2) \\
    = g(\nabla^g u,e_1)\,\mathrm{Hess}\,u(e_1,e_2)
       + g(\nabla^g u,e_2)\,\bigl(-\Delta_g u - \mathrm{Hess}\,u(e_1,e_1)\bigr) \\
    = (\lambda_1 - k_g)\,g(\nabla^g u,e_1)^2
       - g(\nabla^g u,e_2)\,\bigl(-g(\nabla^g f, \nabla^g u) + \mathrm{Hess}\,u(e_1,e_1)\bigr)\\
    = (\lambda_1 - k_g)\,\bigl(g(\nabla^g u,e_1)^2+g(\nabla^g u,e_2)^2\bigr) + g(\nabla^g u,e_2)\,g(\nabla^g f, \nabla^g u)\\
    = (\lambda_1 - k_g)\,|\nabla^g u|^2 + g(\nabla^g u,e_2)\,g(\nabla^g f, \nabla^g u)
\end{gathered}
\end{equation}
at $p \in \partial M$. Hence inequality \eqref{eq:1.16} implies
\begin{equation}\label{eq:1.17.1.1}
\begin{gathered}
(\lambda_1 - k_g)|\nabla^g u|^2 >- g(\nabla^g u,e_2)\,g(\nabla^g f, \nabla^g u)\\= -g(\nabla^g u,e_2)^2\,g(\nabla^g f, e_2) - g(\nabla^g u,e_1)\,g(\nabla^g u,e_2)\,g(\nabla^g f, e_1)
\end{gathered}
\end{equation}
at $p \in \partial M$. Put $Z = g(\nabla^g u, e_1)$, $Y = g(\nabla^g u, e_2)$,
$a = (\partial f/\partial\vec n)(p)$, and $b = g(\nabla^g f, e_1)$. Then
$|\nabla^g u|^2 = Z^2 + Y^2$, and the right-hand side of~(20) equals
$-(aY^2 + bZY)$. The quadratic form $aY^2 + bZY$ in the variables $(Z, Y)$ has the matrix
\[
Q=\begin{pmatrix}0&b/2\\ b/2&a\end{pmatrix}
\]
with eigenvalues $\tfrac12\big(a\pm\sqrt{a^2+b^2}\big)$. Hence $aY^2+bZY\leq\mu_+(a,b)\,(Z^2+Y^2)$, where $\mu_+(a,b)=\lambda_{\max}(Q)$. Returning to the original notation, we get
\begin{equation}\label{eq:1.17}
(\lambda_1 - k_g)|\nabla^g u|^2 > -\mu_+\!\Big(\dfrac{\partial f}{\partial \vec{n}}(p),\dfrac{\partial f}{\partial \vec{\nu}}(p)\Big)|\nabla^g u|^2.
\end{equation}
Since $p$ is a maximum point of the nonconstant function $|\nabla^gu|^2$, we have $|\nabla^g u|^2(p) > 0$. Dividing \eqref{eq:1.17} by $|\nabla^gu|^2$, we obtain the inequality
\begin{equation*}
\begin{gathered}
\lambda_1 >k_g(p) - \mu_+\!\Big(\dfrac{\partial f}{\partial \vec{n}}(p),\dfrac{\partial f}{\partial \vec{\nu}}(p)\Big) \geq c.
\end{gathered}
\end{equation*}
Note that in Case~\ref{I} we obtain the strict inequality $\lambda_1 > c$ in both subcases \ref{i} and~\ref{ii}.

We now prove the inequality in Case \ref{II}. Since constant eigenfunctions correspond to the zero eigenvalue, $u$ is nonconstant and $|\nabla^gu|^2$ is constant. Therefore $|\nabla^g u|\neq 0$ everywhere. Moreover, since $|\nabla^g u|^2 \equiv \mathrm{const}$, inequality~\eqref{eq:1.0} becomes an equality, which implies
\begin{equation*}
\begin{gathered}
\mathrm{Hess}\,u = 0 \quad \text{and} \quad \operatorname{Ric}_{\sigma}(\nabla^gu, \nabla^gu) = 0 \quad \text{on } M.
\end{gathered}
\end{equation*}

Next, using the condition $L_{\sigma} u = 0$, we obtain
\begin{equation}\label{eq:1.18.1}
\begin{gathered}
g( \nabla^g u, \nabla^g f ) = L_{\sigma} u-\Delta_g u =-\Delta_g u= \operatorname{tr}(\operatorname{Hess}u) = \operatorname{tr}(0)=0.
\end{gathered}
\end{equation}
Furthermore,
\begin{equation*}
\begin{gathered}
0 = \operatorname{Hess}u(\nabla^g u, X) = g(\nabla_{\nabla^g u}\nabla^g u, X) \quad \text{at every point } x \in M \text{ and for every } X \in T_xM,
\end{gathered}
\end{equation*}
which gives $\nabla_{\nabla^g u}\nabla^g u=0$. Then
\begin{equation*}
\begin{gathered}
\operatorname{Hess} f(\nabla^gu, \nabla^gu) = g(\nabla_{\nabla^g u}\nabla^g f, \nabla^g u) = \nabla^g u\bigl(g(\nabla^g u, \nabla^g f)\bigr) - g(\nabla^g f, \nabla_{\nabla^g u}\nabla^g u) = 0.
\end{gathered}
\end{equation*}
Consequently, the Gaussian curvature $K$ of the surface satisfies
\begin{equation*}
\begin{gathered}
K\, g( \nabla^g u, \nabla^g u ) = \operatorname{Ric}(\nabla^gu, \nabla^gu) = \operatorname{Ric}_{\sigma}(\nabla^gu, \nabla^gu) = 0,
\end{gathered}
\end{equation*}
and therefore $K \equiv 0$.

Using Proposition~\ref{boundary_cor}, we choose a local orthonormal frame $\{e_1, e_2\}$ along $\partial M$ such that $e_1$ is tangent to $\partial M$ and $e_2 = \vec{n}$ is the outward normal. Then
\begin{equation}\label{eq:1.19}
\begin{gathered}
0 = \operatorname{Hess}u(e_1, e_2)
= e_1 (e_2u) - \nabla_{e_1} e_2(u)
= e_1(\lambda_1 u) - k_g e_1(u)
= (\lambda_1 - k_g)\,e_1(u),
\end{gathered}
\end{equation}
where we used the equality \eqref{eq:1.6}. It follows that $\lambda_1 = k_g$ at every point where $e_1(u) \neq 0$.

We also have
\begin{equation}\label{eq:1.20}
\begin{gathered}
0 = \operatorname{Hess}u(e_1, e_1)
= e_1 (e_1 u) - \nabla_{e_1} e_1(u)
= e_1(e_1 u) + k_g e_2(u)
= e_1(e_1 u) + k_g\, \lambda_1 u,
\end{gathered}
\end{equation}
where we used the equality \eqref{eq:1.5} and the fact that $u$ is an eigenfunction of the Steklov problem~\eqref{eq:1.1}.

\begin{lemma}\label{e_1_equiv_lemma}
The boundary $\partial M$ is connected.
\end{lemma}
\begin{proof}
We have already proved that $K \equiv 0$ on~$M$. Moreover, the assumption \eqref{eq:1.4} implies that $k_g \geq c > 0$ on $\partial M$, because $\mu_+\geq0$. By the Gauss--Bonnet theorem, we get the equality
\[
2\pi\,\chi(M) = \int_M K\,dv_g + \int_{\partial M} k_g\,ds = \int_{\partial M} k_g\,ds > 0.
\]Hence $\chi(M) \geq 1$.

On the other hand, $M$ is a compact connected surface with nonempty boundary. Suppose first that $M$ is nonorientable. Then $M$ is homeomorphic to a connected sum of $k \geq 1$ copies of the projective plane with $b \geq 1$ open disks removed, and $\chi(M) = 2 - k - b$ (see, e.g., \cite[Chapter~1]{Massey}). Hence $\chi(M) \leq 2 - 1 - 1 = 0$, contradicting $\chi(M) \geq 1$. Thus $M$ is orientable and $\chi(M)=2-2\tilde{\gamma}-b$, where $\tilde{\gamma}\geq 0$ is the genus of the surface and $b\geq1$ is the number of connected components of $\partial M$ (see, e.g., \cite[Section~1.1]{farb}). Hence $\chi(M) \leq 1$.

Then the classification theorem for surfaces (see, e.g., \cite[Theorem~1.1]{farb}) implies that $M$ is homeomorphic to a closed two-dimensional disk and $\partial M$ is connected.
\end{proof}

\begin{lemma}\label{e_1_lemma}
The equality $e_1(u) = 0$ holds at only finitely many points of $\partial M$.
\end{lemma}

\begin{proof}
By Lemma~\ref{e_1_equiv_lemma}, $\partial M$ is connected. Let $\hat{u} := u|_{\partial M}$ be the restriction of $u$ to $\partial M$, written in terms of the arc-length parameter $t \in [0, \ell]$ along $\gamma_{\partial M}$ (see Proposition~\ref{boundary_cor}). By \eqref{eq:1.20}, $\hat{u}$ satisfies the ODE
\begin{equation*}
        \ddot{\hat{u}}(t) + k_g\!\left(\gamma_{\partial M}(t)\right)\lambda_1\hat{u}(t) = 0,
        \qquad \hat{u}(0) = \hat{u}(\ell),\quad \dot{\hat{u}}(0) = \dot{\hat{u}}(\ell).
\end{equation*}
Moreover, the assumption~\eqref{eq:1.4} implies that $k_g(x) > 0$ for every point $x \in \partial M$, because $\mu_+\geq0$.
Suppose that there are infinitely many points where $e_1(u) = 0$. By compactness of the interval, there exists a strictly monotone sequence $[0, \ell] \ni t_n \to t_0$ such that $\dot{\hat{u}}(t_n) = 0$. Hence $\dot{\hat{u}}(t_0) = 0$, and by Rolle's theorem there exists a sequence $[0, \ell] \ni t_n' \to t_0$ such that $\ddot{\hat{u}}(t_n') = 0$. Therefore $\ddot{\hat{u}}(t_0) = 0$ and
\[
        \hat{u}(t_0) = -\frac{\ddot{\hat{u}}(t_0)}{k_g\!\left(\gamma_{\partial M}(t_0)\right)\lambda_1} = 0.
\]
Then, by the existence and uniqueness theorem for the Cauchy problem, we conclude that $\hat{u} = u|_{\partial M} \equiv 0$ on $\partial M$. Since $L_{\sigma}u = 0$, Proposition~\ref{strong_max} (applied to $u$ and $-u$) shows that both the minimum and the maximum of $u$ are attained on $\partial M$. Therefore $u \equiv 0$ on $M$, which contradicts the fact that $u$ is a nonzero eigenfunction corresponding to $\lambda_1$.
\end{proof}

By Lemma~\ref{e_1_lemma} and formula~\eqref{eq:1.19}, the equality $\lambda_1 = k_g$ holds on $\partial M$ except at finitely many points. By continuity of $k_g$, we conclude that $\lambda_1 = k_g(x)$ at every point $x \in \partial M$.

Thus, from the assumption \eqref{eq:1.4}, we obtain the inequality
\begin{equation*}\label{eq:1.21}
\begin{gathered}
\lambda_1 = k_g(x) = k_g(x) -\mu_+\!\left(\dfrac{\partial f}{\partial \vec{n}}(x),\dfrac{\partial f}{\partial \vec{\nu}}(x)\right)  + \mu_+\!\left(\dfrac{\partial f}{\partial \vec{n}}(x),\dfrac{\partial f}{\partial \vec{\nu}}(x)\right) \geq c ,
\end{gathered}
\end{equation*}
where the equality $\lambda_1 = c$ holds if and only if $k_g = c$ and
\begin{equation}\label{eq:1.21.1}
\mu_+\!\left(\dfrac{\partial f}{\partial \vec{n}}(x),\dfrac{\partial f}{\partial \vec{\nu}}(x)\right)  = 0 \quad \text{at every point } x \in \partial M.
\end{equation}

It remains to prove the following lemma.

\begin{lemma}\label{last_lemma}
Let $\bar{B}(0, c^{-1})$ denote the closed Euclidean disk of radius $c^{-1}$ centered at the origin. If $\lambda_1 = c$, then $M$ is isometric to $\bar{B}(0, c^{-1})$ and $f$ is constant on $M$.
\end{lemma}
\begin{proof}
Since $\lambda_1=c$, Case \ref{I} is impossible (in Case~\ref{I} we proved $\lambda_1>c$). Hence equality~\eqref{eq:1.18.1} holds and problem~\eqref{eq:1.1} for the eigenfunction $u$ reduces to the problem
\begin{equation*}\label{eq:1.22}
\begin{cases}
\Delta_g u = 0 & \text{in } M, \\
\dfrac{\partial u}{\partial \vec{n}} = \lambda_1 u & \text{on } \partial M.
\end{cases}
\end{equation*}
Moreover, the equality $\lambda_1=c$ implies that $k_g = c  > 0$ on $\partial M$. In addition, $K\equiv0$. Hence, by Theorem~\ref{escobar_theor}, there is an isometry $\Phi\colon M\to\bar{B}(0,c^{-1})$. Identifying $M$ with $\bar{B}(0,c^{-1})$ via $\Phi$ (and replacing $f$, $u$ by $f\circ\Phi^{-1}$, $u\circ\Phi^{-1}$), we may assume that $M$ is the Euclidean disk with polar coordinates $(r,\psi)$.

Thus $u$ is an eigenfunction of the classical Steklov problem on the two-dimensional Euclidean disk of radius $c^{-1}$ corresponding to the eigenvalue $c$. Then Examples 1.3.1 and 1.3.2 in~\cite{girouard} give $u = r(\alpha \cos\psi+\beta\sin \psi)$ in polar coordinates for some pair of constants $(\alpha,\beta) \neq (0,0)$.

We now find all functions $f$ satisfying conditions~\eqref{eq:1.18.1} and~\eqref{eq:1.21.1}. Put $R=c^{-1}$, $x=r\cos\psi$, and $y=r\sin\psi$, so that $u=r(\alpha\cos\psi+\beta\sin\psi)=\alpha x+\beta y$. Let $A=\sqrt{\alpha^2+\beta^2}>0$ and pass to the new coordinates
\[
x'=\frac{\alpha x+\beta y}{A},\qquad y'=\frac{-\beta x+\alpha y}{A},
\]
that is, rotate the plane by the angle $\theta$ with $\cos\theta=\alpha/A$, $\sin\theta=\beta/A$. In these coordinates $u=Ax'$. The rotation is an isometry of the disk. It maps the boundary circle onto itself and the outward normal onto the outward normal, so the normal and tangential derivatives of $f$ (the latter up to sign) on $\partial M$ are the same in the new coordinates. We keep writing $(x,y)$ instead of $(x',y')$,
so $u=Ax$.

Condition~\eqref{eq:1.18.1}, that is $g(\nabla^g u,\nabla^g f)=0$, gives $A\,\partial_x f=0$ on $M$, so $\partial_x f\equiv0$. Each horizontal section of the disk is an interval, therefore $f(x,y)=\tilde{f}(y)$ with $\tilde{f}(y)$ smooth on $[-R,R]$ and $\tilde{f}(y)=f(0,y)$.

By the definition of $\mu_+$, condition~\eqref{eq:1.21.1} gives $\dfrac{\partial f}{\partial \vec\nu}=0$ on $\partial M$. On the circle $x^2+y^2=R^2$, the unit tangent vector is $\vec\nu=\pm(-y,x)/R$, so
\[
0=\dfrac{\partial f}{\partial \vec\nu}=\pm\frac{x}{R}\,\tilde{f}'(y).
\]
Hence $\tilde{f}'(y)=0$ for every $y\in(-R,R)$, since each such $y$ is the ordinate of a point of $\partial M$ with $x\neq0$. Therefore $\tilde{f}$ is constant on $(-R,R)$, and, by continuity, $\tilde{f}$ is constant on $[-R,R]$. So $f$ is constant on $M$.
\end{proof}

It remains to prove the converse implication. Suppose that $M$ is isometric to $\bar{B}(0, c^{-1})$ and $f \equiv \mathrm{const}$. Then $L_{\sigma} = \Delta_g$, and Examples 1.3.1 and 1.3.2 in~\cite{girouard} give $\lambda_1 = c$. This completes the proof of Theorem \ref{main_t}.

\bibliographystyle{plain}
\bibliography{references}

\end{document}